\documentclass{article}

\usepackage{amssymb}

\usepackage{graphicx}
\usepackage{amscd, color}

\newtheorem{defi}{Definition}[section]
\newtheorem{theo}{Theorem}[section]

\newtheorem{remark}[theo]{Remark}

\newtheorem{lemma}[theo]{Lemma}

\newtheorem{coro}[theo]{Corollary}
\newtheorem{con}[theo]{Conjecture}
\newtheorem{prop}[theo]{Proposition}
\newtheorem{fact}[theo]{Fact}
\newtheorem{algorithm}[theo]{Algorithm}

\newcommand{\qed}{\hspace*{\fill} \rule{7pt}{7pt}}

\begin{document}

\title{On Cliques and Lagrangians of 3-uniform Hypergraphs\footnote{Supported by National Natural Science Foundation of China (No. 11271116)}}

\date{}

\author{ Yuejian Peng \thanks{ College of Mathematics, Hunan University, Changsha 410082, P.R. China.  Email: ypeng1@163.com. } \and Hegui Zhu \thanks{College of Sciences, Northeastern University, Shenyang, 110819, China, and School of Mathematics, Jilin University, Changchun, 130012, China. Email: zhuhegui@yahoo.com.cn}\and Yanling Zheng\thanks{ College of Mathematics, Hunan University, Changsha 410082, P.R. China. }
\and Cheng Zhao \thanks{Department of Mathematics and Computer Science, Indiana State University, Terre Haute, IN, 47809 and School of Mathematics, Jilin University, Changchun 130022, P.R. China. Email: cheng.zhao@indstate.edu }
}


\maketitle
\date

\begin{abstract}
A remarkable connection between the maximum clique number  and the Lagrangian of a graph  was established by Motzkin and Straus. This connection and its extensions were sucessfully employed in optimization to provide heuristics for the maximum clique number in graphs.  In this paper, we provide evidence that the Lagrangian of a 3-uniform hypergraph is related to the order of its maximum cliques when the number of edges of the hypergraph is in certain ranges. In particular,  we present some results about a  conjecture related to Frankl-F\"uredi's conjecture about Lagrangians of hypergraphs. We also describe a combinatorial algorithm that can be used to check the validity of the conjecture.
\end{abstract}

keywords: Cliques of Hypergraphs; Lagrangians of Hypergraphs; Optimization.

\section{Introduction}
In 1965,  Motzkin and Straus \cite{MS} provided a new proof of Tur\'an's theorem based on a continuous characterization of the clique number of a graph using the Lagrangian of a graph. This new proof aroused interests in the study of Lagrangians of hypergraphs. Furthermore, the Motzkin-Straus result and its extension were successfully employed in optimization to provide heuristics for the maximum clique problem, and  the Motzkin-Straus theorem has  been also generalized to vertex-weighted graphs \cite{G9} and edge-weighted graphs with applications to pattern recognition in image analysis (see \cite{B1}, \cite{B2}, \cite{B3}, \cite{G9}, \cite{PP}, \cite{PP15}, \cite{RTP20}).   
In this paper, we provide evidence that the Lagrangian of an $r$-uniform hypergraph is related to the order of its maximum cliques under some conditions. We first state a few definitions.

For a set $V$ and a positive integer $r$ we denote by $V^{(r)}$ the family of all $r$-subsets of $V$. An $r$-uniform graph or $r$-graph $G$ consists of a set $V(G)$ of vertices and a set $E(G) \subseteq V(G) ^{(r)}$ of edges. An edge $e=\{a_1, a_2, \ldots, a_r\}$ will be simply denoted by $a_1a_2 \ldots a_r$. The complement of an $r$-graph $G$ is denoted by $G^c$.  An $r$-graph $H$ is  a {\it subgraph} of an $r$-graph $G$, denoted by $H\subseteq G$ if $V(H)\subseteq V(G)$ and $E(H)\subseteq E(G)$.  Let ${\mathbb N}$ be the set of all positive integers. For any integer $n \in {\mathbb N}$,  denote the set $\{1, 2, 3, \ldots, n\}$ by $[n]$. Let $K^{(r)}_t$ denote the complete $r$-graph on $t$ vertices, that is the $r$-graph on $t$ vertices containing all possible edges. A complete $r$-graph on $t$ vertices is also called a clique with order $t$. We also let $[n]^{(r)}$  represent the  complete $r$-uniform graph on the vertex set $[n]$. When $r=2$, an $r$-uniform graph is a simple graph.  When $r\ge 3$,  an $r$-graph is often called a hypergraph.


\begin{defi}
Let $G$ be  an $r$-uniform graph  with vertex set $\{1,2,\ldots,n\}$ and
edge set $E(G)$.  Let $S=\{\vec{x}=(x_1,x_2,\ldots ,x_n)\in R^n: \sum_{i=1}^{n} x_i =1, x_i
\ge 0 {\rm \ for \ } i=1,2,\ldots , n \}$. For $\vec{x}=(x_1,x_2,\ldots ,x_n)\in S$,
define
$$\lambda (G,\vec{x})=\sum_{i_1i_2 \cdots i_r \in E(G)}x_{i_1}x_{i_2}\ldots x_{i_r}.$$
The Lagrangian of
$G$, denoted by $\lambda (G)$, is defined as $\lambda (G) = \max \{\lambda (G, \vec{x}): \vec{x} \in S \}$.
 A vector $\vec{y}\in S$ is called an {\em optimal weighting} for $G$ if $\lambda (G, \vec{y})=\lambda(G)$.
\end{defi}

 The following fact is easily implied by the definition of the Lagrangian.

\begin{fact}\label{mono}
Let $G_1$, $G_2$ be $r$-uniform graphs and $G_1\subseteq G_2$. Then $\lambda (G_1) \le \lambda (G_2).$
\end{fact}

In \cite{MS}, Motzkin and Straus proved that  the Lagrangian of a 2-graph is determined by the order of its maximum clique.

\begin{theo} (Motzkin and Straus \cite{MS}) \label{MStheo}
If $G$ is a 2-graph in which a largest clique has order $l$ then
$\lambda(G)=\lambda(K^{(2)}_l)=\lambda([l]^{(2)})={1 \over 2}(1 - {1 \over l})$.
\end{theo}

The obvious generalization of Motzkin and Straus' result to hypergraphs is false because there are many examples of hypergraphs that do not achieve their Lagrangian on any proper subhypergraph.
An attempt to generalize the Motzkin-Straus theorem to hypergraphs is due to S\'os and Straus \cite{SS}.  Recently, in \cite{BP1} and \cite{BP2} Rota Bul\'o and Pelillo generalized the Motzkin and Straus' result to $r$-graphs in some way using a continuous characterization of maximal cliques with applications in image analysis. Although, the obvious generalization of Motzkin and Straus' result to hypergraphs is false,  we attempt to explore the relationship between the Lagrangian of a hypergraph and the order of its maximum cliques for hypergraphs when the number of edges is in certain range. In \cite{PZ}, the following two conjectures are proposed.

\begin{con} \label{conjecture1} (Peng-Zhao \cite{PZ})
Let $l$ and $m$ be positive integers satisfying ${l-1 \choose r} \le m \le {l-1 \choose r} + {l-2 \choose r-1}$.
Let $G$ be an $r$-graph with $m$ edges containing a clique of order  $l-1$. Then $\lambda(G)=\lambda([l-1]^{(r)})$.
\end{con}

\begin{con} \label{conjecture2} (Peng-Zhao \cite{PZ})
 Let $G$ be an $r$-graph with $m$ edges and without containing a clique of order  $l-1$, where ${l-1 \choose r} \le m \le {l-1 \choose r} + {l-2 \choose r-1}$. Then $\lambda(G) < \lambda([l-1]^{(r)})$.
\end{con}

The upper bound ${l-1 \choose r} + {l-2 \choose r-1}$ in Conjecture \ref{conjecture1} is the best possible. For example, if $m ={l-1 \choose r}+{l-2 \choose r-1}+1$ then $\lambda(C_{r,m}) > \lambda([l-1]^{(r)})$, where  $C_{r,m}$ is the $r$-graph on the  vertex set $[l]$ and with the edge set $[l-1]^{(r)}\cup\{i_1\cdots i_{r-1}l, i_1\cdots i_{r-1}\in [l-2]^{(r-1)}\}\cup\{1\cdots (r-2)(l-1)l\}$.

In the course of estimating Tur\'an densities of hypergraphs by applying the Lagrangians of related hypergraphs, Frankl and F\"uredi \cite{FF} asked the following question: Given $r \ge 3$ and $m \in {\mathbb N}$ how large can the Lagrangian of an $r$-graph with $m$ edges be? In order to state their conjecture on this problem we require the following definition. For distinct $A, B \in {\mathbb N}^{(r)}$ we say that $A$ is less than $B$ in the {\em colex ordering} if $max(A \triangle B) \in B$, where $A \triangle B=(A \setminus B)\cup (B \setminus A)$ is the symmetric difference of $A$ and $B$. For example we have $246 < 156$ in ${\mathbb N}^{(3)}$ since $max(\{2,4,6\} \triangle \{1,5,6\}) \in \{1,5,6\}$.  Let $C_{r,m}$ denote the $r$-graph with $m$ edges formed by taking the first $m$ elements in the colex ordering of ${\mathbb N}^{(r)}$. The following conjecture of Frankl and F\"uredi (if it is true) proposes a  solution to the above question.

\begin{con} (Frankl and F\"uredi \cite{FF})\label{conjecture} The $r$-graph with $m$ edges formed by taking the first $m$ sets in the colex ordering of ${\mathbb N}^{(r)}$ has the largest Lagrangian of all $r$-graphs with $m$ edges. In particular, the $r$-graph with $l \choose r$ edges and the largest Lagrangian is $[l]^{(r)}$.
\end{con}

This conjecture is true when $r=2$ by Theorem \ref{MStheo}. For the case $r=3$, Talbot in \cite{T} proved the following.

\begin{theo} (Talbot \cite{T}) \label{Tal} Let $m$ and $l$ be integers satisfying
$${l-1 \choose 3} \le m \le {l-1 \choose 3} + {l-2 \choose 2} - (l-1).$$
Then Conjecture \ref{conjecture} is true for $r=3$ and this value of $m$.  Conjecture \ref{conjecture} is also true for $r=3$ and $m= {l \choose 3}-1$ or $m={l \choose 3} -2$.
\end{theo}

The truth of Frankl and F\"uredi's conjecture is not known in general for $r \ge 4$. Even in the case $r=3$, Theorem \ref{Tal} does not cover the case when ${l-1 \choose 3}+{l-2 \choose 2}-(l-2) \le m \le {l \choose 3}-3$ in this conjecture. In \cite{HPZ}, He, Peng, and Zhao verified Frankl and F\"uredi's conjecture  for some values $m$ when $r=3$.

 The following result is  given in \cite{T}.

\begin{lemma} (Talbot \cite{T}) \label{LemmaTal7}
For any integers $m,l,$ and $r$ satisfying ${l-1 \choose r} \le m \le {l-1 \choose r} + {l-2 \choose r-1}$,
we have $\lambda(C_{r,m}) = \lambda([l-1]^{(r)})$.
\end{lemma}

If Conjectures \ref{conjecture1} and \ref{conjecture2} are true, then Conjecture \ref{conjecture} is true for this range of $m$.
In \cite{PZ}, it has been shown that  Conjecture \ref{conjecture1} holds when $r=3$.

\begin{theo} \label{theorem 1} (Peng-Zhao \cite{PZ}) Let $m$ and $l$ be positive integers satisfying ${l-1 \choose 3} \le m \le {l-1 \choose 3} + {l-2 \choose 2}$. Let $G$ be a $3$-graph with $m$ edges and $G$ contain a clique of order  $l-1$. Then $\lambda(G) = \lambda([l-1]^{(3)})$.
\end{theo}

Further evidences for Conjectures \ref{conjecture1} and \ref{conjecture2} are provided  in \cite{PTZ} when the number of vertices of the r-graph is restricted to be $l$.

\begin{defi}
An $r$-graph $G=([n],E)$ is {\it left-compressed} if $j_1j_2 \cdots j_r \in E$ implies $i_1i_2 \cdots i_r \in E$ provided $i_p \le j_p$ for every $p, 1\le p\le r$. 
\end{defi}
\textcolor{blue}{In \cite{T}, Talbot showed that to confirm Conjecture \ref{conjecture}, it is sufficient to verify for left-compressed $r$-graphs. (Lemma 2.3 in \cite{T}). The proof of this reuction lemma is  to start with an $r$-graph with $m$ edges which has the largest Lagrangian of all $r$-graphs with $m$ edges.Such an $r$-graph is called an extremal $r$-graph  for $m$. Then perform 
a sequence of `left-compressing' operations  to $G$ (for $j_1j_2 \cdots j_r \in E$, if there exists $i_1i_2 \cdots i_r \notin E$, where $i_p \le j_p$ for every $p$, then replace $j_1j_2 \cdots j_r$ by $i_1i_2 \cdots i_r$ until there is no such an edge) to get a left-compressed extremal $r$-graph for $m$. This reduction lemma is very useful and has been applied in several later work (need to add references here).
In this paper, we  first show a similar reduction lemma for $r=3$: we only need to consider left-compressed  3-graphs  to confirm Conjecture \ref{conjecture2} when $r=3$ (Theorem \ref{theorem2}). However, the argument in \cite{T} doesn't trivally apply to the reduction lemma in this situation.  If we start with a $3$-graph with with $m$ edges and without containing a clique of order  $l-1$, the resulting left-compressed $3$-graph obtained by a sequence of left-compressing operations might contain a clique of order  $l-1$. So the proof of the reduction lemma in this situation is not as easy as we thought in the beginning. In fact,  a preliminary related partial results was given in \cite{PZa}. However the proof of Lemma 3 in  \cite{PZa} (Lemma \ref{lemma1} in this paper, the key ingredient in the proof of Theorem \ref{theorem2}) was  flawed. Theorem \ref{theorem2} allows us to check  only for left-compressed 3-graphs on vertex set $[l]$ to verify  Conjecture \ref{conjecture2} (for $r=3$) and it builds a foundation for  some arguments in other papers (\cite{SPT} and \cite{TPZZ2}).  So it is important to have a correct proof. }

To make it accurate, we have

\begin{theo}\label{theorem2} To verify Conjecture \ref{conjecture2} for $r=3$ and  any given $l$, it is sufficient to verify $\lambda (G)<\lambda ([l-1]^{(3)})$  for all left-compressed $3$-graphs $G$ on the vertex set $[l]$  with $m = {l-1 \choose 3} + {l-2 \choose 2}$ edges and without containing the clique $[l-1]^{(3)}$.
\end{theo}

Based on Theorem \ref{theorem2}, we  describe an algorithm of verifying Conjecture \ref{conjecture2} for  given $l$ in Section \ref{prooftheorem2}. 

As an implication, we also show that

\begin{coro} \label{mbc}
Let $G$ be a $3$-graph  with $m$ edges and  containing no clique of order  $l-1$,
where ${l-1 \choose 3}\le m \le {l-1 \choose 3} + {l-2 \choose 2}$.  If $6 \le l \le 13$, then $\lambda(G) < \lambda([l-1]^{(3)})$.
\end{coro}

Corollary \ref{mbc} is applied to prove a result in \cite{TPZZ4}.

The following result provide more evidence for  Conjecture \ref{conjecture2} when $r=3$.

\begin{theo} \label{Lemma103} Let $G$  be  a left-compressed $3$-graph on $[l]$ with $m = {l-1\choose 3}+{l-2\choose 2}$ edges.  If the first $j+1$ triples in colex ordering in $G^{c}=[l]^3 \backslash E(G)$ are  $(l-2-j)(l-2)(l-1), (l-2-j+1)(l-2)(l-1), \cdots, (l-3)(l-2)(l-1)$ and $ (l-2-i)(l-2)l$, where $i\ge j\ge 1$,   then $ \lambda(G)<\lambda ([l-1]^{(3)})$.
\end{theo}

The proofs of Theorem \ref{theorem2}, Corollary \ref{mbc}, and Theorem \ref{Lemma103} will be given in Sections \ref{prooftheorem2}, \ref{proofmbc}, and \ref{proofLemma103} respectively.

Let us state some preliminary results in the following section.

\section{Preliminary Results}
For an $r$-graph $G=(V,E)$ we denote the $(r-1)$-neighborhood of a vertex $i \in V$ by $E_i=\{A \in V^{(r-1)}: A \cup \{i\} \in E\}$. Similarly, we will denote the $(r-2)$-neighborhood of a pair of vertices $i,j \in V$ by $E_{ij}=\{B \in V^{(r-2)}: B \cup \{i,j\} \in E\}$. We denote the complement of $E_i$ by $E^c_i=\{A \in V^{(r-1)}: A \cup \{i\} \in V^{(r)} \backslash E\}$. Also, we will denote the complement of $E_{ij}$ by
$E^c_{ij}=\{B \in V^{(r-2)}: B \cup \{i,j\} \in V^{(r)} \backslash E\}$. Denote $$E_{i\setminus j}=E_i\cap E^c_j.$$   In some case, we will impose one additional condition on any optimal weighting ${\vec x}=(x_1, x_2, \ldots, x_n)$ for an $r$-graph $G$:
\begin{equation}
  |\{i : x_i > 0 \}|{\rm \ is \ minimal, i.e. \ if}  \ \vec y \in S
  {\rm \ satisfying } \ \  |\{i : y_i > 0 \}| < |\{i : x_i > 0 \}|,
 {\rm \  then \ } \lambda (G, {\vec y}) < \lambda(G) \label{conditionb}.
\end{equation}

When the theory of Lagrange multipliers 
is applied to find the optimum of $\lambda(G)$, subject to $\sum_{i=1}^n x_i =1$, notice that $\lambda (E_i, {\vec x})$ corresponds to the partial derivative of  $\lambda(G, \vec x)$ with respect to $x_i$.
The following lemma gives some necessary condition of an optimal weighting of  $\lambda(G)$.

\begin{lemma} (Frankl and R\"{o}dl \cite{FR84}) \label{LemmaTal5} \textcolor{blue}{Let $G=(V,E)$ be an $r$-graph on the vertex set $[n]$ and ${\vec x}=(x_1, x_2, \ldots, x_n)$ be an optimal weighting for $G$ with $k$  ($\le n$) positive weights $x_1, x_2, \cdots, x_k$. Then for every $\{i, j\} \in [k]^{(2)}$,}

(a) $\lambda (E_i, {\vec x})=\lambda (E_j, \vec{x})=r\lambda(G)$.

\textcolor{blue}{(b) If ${\vec x}$ satisfies condition (\ref{conditionb}), then there is an edge in $E$ containing both $i$ and $j$.}
\end{lemma}

\begin{lemma}\label{contain}
Suppose $G$ is an $r$-uniform graph on the vertex set  $[n]$ with edge set $E$. Let $1\le i<j\le n$.
If  $E_{j\setminus i}=\emptyset$, then there exists an optimal weighting $\vec{y}=(y_1, y_2, \ldots, y_n)$  of $\lambda (G)$ such that  $y_i\ge y_j$.
\end{lemma}

{\em Proof.} Let  $\vec{x}=(x_1, x_2, \ldots, x_n)$ be an optimal weighting for $G$. If $x_i<x_j$, then let
$y_k=x_k$ for $k\neq i, j$, $y_i=x_j$ and $y_j=x_i$. Then $\vec{y}=(y_1, y_2, \ldots, y_n)$ is a  weighting for $G$ with $y_i>y_j$ and
$$\lambda(G, \vec y)-\lambda(G, \vec x)=(x_j-x_i)(\lambda (E_{i\setminus j}, {\vec x})-\lambda (E_{j\setminus i}, {\vec x}))\ge 0.$$
So $\vec{y}$ is an optimal weighting satisfying the condition.
\qed

We call two vertices $i,j$ of an $r$-uniform graph $G=(V, E)$ equivalent if for all $f \in {V-\{i,j\} \choose r-1}$, $f \in E_i$ if and only if $f  \in E_j$.


\begin{lemma}\label{forlambda}(c.f. \cite{FR84})
Suppose $G$ is an $r$-uniform graph on the vertex set $[n]$.
If vertices $i_1$, $i_2$, ..., $i_t$
 are pairwisely equivalent,  then there exists an optimal weighting $\vec{y}=(y_1, y_2, \ldots, y_n)$  of $\lambda (G)$ such that  $y_{i_1} = y_{i_2}=\cdots =y_{i_t}$.
\end{lemma}

We note that an $r$-graph $G=(V,E)$ on the  vertex set $[n]$ is left compressed if and only if $E_{j\setminus i}=\emptyset$ for any $1\le i<j\le n$.

\begin{remark}\label{r1} \textcolor{blue}{ Let $G=(V,E)$ be an $r$-graph on the vertex set $[n]$ and ${\vec x}=(x_1, x_2, \ldots, x_n)$ be an optimal weighting for $G$ with $k$  ($\le n$) positive weights $x_1, x_2, \cdots, x_k$. Let $1\le i<j\le k$.  Then}

(a) In Lemma \ref{LemmaTal5}, part(a) implies that
$x_j\lambda(E_{ij}, {\vec x})+\lambda (E_{i\setminus j}, {\vec x})=x_i\lambda(E_{ij}, {\vec x})+\lambda (E_{j\setminus i}, {\vec x})$.
\textcolor{blue}{In particular, if $E_{j\setminus i}=\emptyset$,} then
\begin{equation}\label{enbhd}
(x_i-x_j)\lambda(E_{ij}, {\vec x})=\lambda (E_{i\setminus j}, {\vec x}).
\end{equation}

(b) By (\ref{enbhd}), if  $G$ is left-compressed, then  ${\vec x}$ must satisfy
\begin{equation}\label{conditiona}
x_1 \ge x_2 \ge \ldots \ge x_n \ge 0.
\end{equation}
\end{remark}

\section{Proof of Theorem \ref{theorem2}}\label{prooftheorem2}

Denote $\lambda_{(m,l)}^{r-}=\max \{ \lambda(G): G \ is {\rm \ an \ } r-{\rm graph \ with \ } m   {\rm \ edges } {\rm \ and \ } {\rm \ does \ not \ contain \ a \ clique \ of \ size \ }   l  \}$ .

We need the following reduction lemma.

\begin{lemma} \label{lemma1} Let $m$ and $l$  be  positive integers satisfying ${l-1 \choose 3} \le m \le {l-1 \choose 3} + {l-2 \choose 2}$.
Then there exists a left compressed $3$-graph $G$ on the vertex set $[l]$ with $m$ edges  such that $\lambda(G)=\lambda_{(m,l-1)}^{3-}$ or Conjecture \ref{conjecture2} holds.
\end{lemma}

To verify this lemma, we define a partial order Hessian diagram $K$ (2-graph on vertices of all possible triples $i_1i_2i_3$ where $1 \le i_1<i_2 <i_3 \le l$). A triple $i_1 i_2i_3$ is called an {\it ancestor } of a triple  $j_1j_2j_3$ if $i_1\ge j_1$, $i_2\ge j_2$, $i_3 \ge j_3$, and $i_1+i_2+i_3 > j_1+j_2+j_3$. In this case, the triple $j_1j_2j_3$  is called a {\it descendant} of $i_1i_2i_3$.  We say that $i_1i_2i_3$ has higher hierarchy than $j_1j_2j_3$ if $i_1i_2i_3$ is  an  ancestor of  $j_1j_2j_3$.  A triple $i_1 i_2i_3$ is called a  {\it directed ancestor } of a triple  $j_1j_2j_3$ if $i_1\ge j_1$, $i_2\ge j_2$, $i_3 \ge j_3$, and $i_1+i_2+i_3 = j_1+j_2+j_3+1$. In this case, the triple $j_1j_2j_3$ is called a {\it directed descendant} of $i_1i_2i_3$, and vertices $i_1 i_2i_3$ and $j_1j_2j_3$ are adjacent in the corresponding Hessian diagram. Note that hierarchy is a partial ordering. Figure 1 shows  part of the hierarchy relationship of triples in $[k]^{(3)}$. Note that a $3$-graph $G$ is left-compressed if and only if for any edge in $G$, all its descendants (in Hessian diagram $K$) should be in $G$ as well.

\begin{figure}[!t]
\centering
\includegraphics[width=4in]{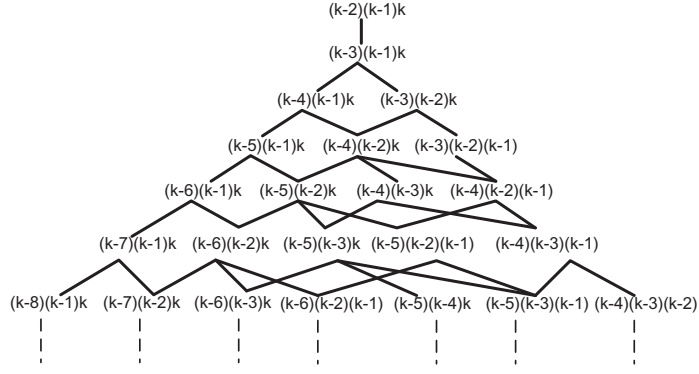}
\caption{Partial order Hessian diagram}
\label{fig_sim}
\end{figure}

\bigskip

\noindent {\bf Proof of Lemma \ref{lemma1}.} Let $G$ be a 3-graph on the vertex set $[n]$ with $m$ edges without containing a clique of order $l-1$  such that $\lambda(G)=\lambda_{(m,l-1)}^{3-}$. We call such a $G$ an extremal 3-graph for $m$ and $l-1$. Let ${\vec x}=(x_1, x_2, \ldots, x_n)$ be an optimal weighting of $G$ and and $k$ be the number of non-zero weights
in $\vec{x}$. If $k\leq l-1$,  Conjecture holds since $G$ does not contain a clique order of $l-1$. So we assume $k\geq l.$  We can assume that $x_i\ge x_j$ when $i<j$ since otherwise we can just relabel the vertices of $G$ and obtain another  extremal $3$-graph for $m$ and $l-1$ with an optimal weighting ${\vec x}=(x_1, x_2, \ldots, x_n)$ satisfying $x_i\ge x_j$ when $i<j$. Next we obtain a new $3$-graph $H$ from $G$ by performing the following:

\begin{enumerate}

\item If $(l-3)(l-2)(l-1) \in E(G)$, then there is at least one triple  in $[l-1]^{(3)}\setminus E(G)$.  We replace $(l-3)(l-2)(l-1)$ by this triple;

\item If  an  edge in $G$ has  a  descendant other than $(l-3)(l-2)(l-1)$  that is not  in $E(G)$, then replace this edge by a descendant other than $(l-3)(l-2)(l-1)$ with the lowest hierarchy. Repeat this until there is no such an edge.

\end{enumerate}

  Then $H$ satisfies the following properties:

\begin{enumerate}

\item  The number of edges in $H$ is the same as the number of edges in $G$.

\item $\lambda(G)=\lambda(G, {\vec x})\le \lambda(H, {\vec x})\le \lambda(H).$

\item $(l-3)(l-2)(l-1) \notin E(H)$.

\item For any edge in $E(H)$, all its  descendants  other than $(l-3)(l-2)(l-1)$ will be in $E(H)$.

\end{enumerate}

If $H$ is not left-compressed, then there is an ancestor $uvw$ of $(l-3)(l-2)(l-1)$ such that $uvw\in E(H)$. We claim that $uvw$ must be $(l-3)(l-2)l$. If $uvw$ is not $(l-3)(l-2)l$, then since all  descendants  other than $(l-3)(l-2)(l-1)$ of $uvw$ will be in $E(H)$, then all descendants of $(l-3)(l-1)l$ (other than $(l-3)(l-2)(l-1)$) or
all descendants of $(l-3)(l-2)(l+1)$ (other than $(l-3)(l-2)(l-1)$) will be in $E(H)$. So
all triples in $[l-1]^{(3)}\setminus\{(l-3)(l-2)(l-1)\}$, all triples in the form of $ijl$ (where $ij\in [l-2]^{(2)}$), and all triples in the form of $ij(l+1)$ (where $ij\in [l-2]^{(2)}$) or all triples in the form of  $i(l-1)l$, $1\le i\le l-3$ will be in $E(H)$, then
$$m\ge {l-1 \choose 3}-1 + {l-2 \choose 2}+(l-3)>{l-1 \choose 3} + {l-2 \choose 2}$$ which is a contradiction.
So $uvw$ must be $(l-3)(l-2)l$. Since $m\le {l-1 \choose 3} + {l-2 \choose 2}$ and all the descendants other than
$(l-3)(l-2)(l-1)$ of an edge in $H$ will be an edge in $H$,
then there are three possibilities.

 Case 1. $E(H)=([l-1]^{(3)}\setminus \{(l-1)(l-2)(l-3)\})\cup\{ijl, ij\in [l-2]^{(2)}\}\cup\{12(l+1)\}.$

 Case 2.  $E(H)=([l-1]^{(3)}\setminus \{(l-1)(l-2)(l-3)\})\cup\{ijl, ij\in [l-2]^{(2)}\} .$

 Case 3.  $E(H)=([l-1]^{(3)}\setminus \{(l-1)(l-2)(l-3)\})\cup\{ijl, ij\in [l-2]^{(2)}\}\cup \{1(l-1)l\}.$

Let ${\vec y}=(y_1, y_2, \ldots, y_n)$ be an optimal weighting for $H$, where  $n=l+1$ or $n=l$.  We claim that if Case 1 happens, then $y_{l+1}=0$. Notice that $E_{(l+1)\setminus i}=\emptyset$  for each $1\le i\le l$, by Lemma \ref{contain},
we can assume that $y_i\ge y_{l+1}$ for each $1\le i\le l$.
If $y_{l+1}>0$, then each $y_i>0$. This contradicts to
$E_{l(l+1)}=\emptyset$ by Lemma \ref{LemmaTal5}. So we have $y_{l+1}=0$. This implies Case 1 is equivalent to Case 2. So we only need to  consider Case 2 and Case 3. We will show if these two cases happen, then Conjecture \ref{conjecture2} holds.

Note that $([l-1]^{(3)}\setminus \{(l-1)(l-2)(l-3)\})\cup\{ijl, ij\in [l-2]^{(2)}\} \subseteq ([l-1]^{(3)}\setminus \{(l-1)(l-2)(l-3)\})\cup\{ijl, ij\in [l-2]^{(2)}\}\cup \{1(l-1)l\}.$
It is sufficient to assume  $E(H)=([l-1]^{(3)}\setminus \{(l-1)(l-2)(l-3)\})\cup\{ijl, ij\in [l-2]^{(2)}\}\cup \{1(l-1)l\}$ and show  $ \lambda(H,\vec{x})<\lambda([l-1]^{(3)})$ since then $\lambda(G)=\lambda(G,\vec{x})\leq \lambda(H,\vec{x})<\lambda([l-1]^{(3)})$, i.e. Conjecture \ref{conjecture2} holds.

 \textcolor{blue}{Let $H'=[l-1]^{(3)}\cup\{ijl, ij\in [l-2]^{(2)}\}\setminus \{(l-3)(l-2)l\}\cup \{1(l-1)l\}.$ Then $\lambda(H,\vec{x})\leq \lambda(H',\vec{x})$ since $x_1 \ge x_2 \ge \ldots \geq x_{t+1}>0$. Note that $H'$ contains $[l-1]^{(3)}$ and the number of the edges in  $H'$ is ${l-1 \choose 3} + {l-1 \choose 2}$. Hence $\lambda(H')=\lambda([l-1]^{(3)})$ by Theorem \ref{theorem 1}. We claim $\vec{x}$ is not an optimal weight for $H'.$ So $\lambda(H',\vec{x})<\lambda(H')=\lambda([l-1]^{(3)}).$ To show this we prove that an optimal weighting of $H'$ must have  $l-1$ positive weights which contradicts to $\vec{x}$ has $l$ positive weights. Clearly, an optimal weighting for $H'$ has at least   $l-1$ positive weights. Note that $H'$ is left-compressed. Let ${\vec z}=(z_1, z_2, \ldots, z_{l})$ be an optimal weighting of $H'$. Then $z_1 \ge z_2 \ge \ldots \geq z_{l}\geq 0$. Suppose $z_{l}>0$ for a contradiction. Let $H^*=H'\backslash \{1(l-1)l\}\bigcup \{(l-3)(l-2)l\}.$ Since  $H^*$ contains $[l-1]^{(3)}$ and the number of the edges in  $H^*$ is ${l-1 \choose 3} + {l-2 \choose 2},$ we have $\lambda(H^*)=\lambda([l-1]^{(3)})$ by Theorem \ref{theorem 1}.  We will show $\lambda(H',\vec{z})<\lambda(H^*,\vec{z})\leq \lambda([l-1]^{(3)})=\lambda(H')$ when $z_{l}>0.$  This contradicts to $\vec{z}$ is an optimal weighting for $H'$. Hence $z_{l}=0.$ Clearly,}
 \begin{eqnarray}\label{eq421}
 \lambda(H^*,\vec{z})-\lambda(H',\vec{z})=z_{l-3}z_{l-2}z_{l}-z_{1}z_{l-1}z_{l}.
\end{eqnarray}
Using Remark \ref{r1}(a), we have
 \begin{eqnarray}\label{eq422}
 z_1=z_{l-1}+\frac{ (z_2+\ldots+z_{l-2})z_{l}}{z_2+\ldots+z_{l-2}+z_{l}}< z_{l-1}+z_{l},
\end{eqnarray}
 \begin{eqnarray}\label{eq423}
 z_1=z_{l-3}+\frac{ (z_{l-2}+z_{l-1})z_{l}}{z_2+\ldots+z_{l-4}+z_{l-2}+z_{l-1}+z_{l}},
\end{eqnarray}
and
 \begin{eqnarray}\label{eq424}
 z_{l-2}=z_{l-1}+\frac{ (z_{2}+\ldots+z_{l-4})z_{l}}{z_1+\ldots+z_{l-3}}.
\end{eqnarray}

Combing (\ref{eq422}),(\ref{eq423}) and (\ref{eq424}), we have
 \begin{eqnarray}\label{eq425}
0<z_1-z_{l-3}<z_{l-2}-z_{l-1}
\end{eqnarray}
for $l\geq 7$ ( We have $z_{l-3}z_{l-2}z_{l}-z_{1}z_{l-1}z_{l}>0$ for $l\leq 6$ by a direction calculation).
Applying (\ref{eq425}) to (\ref{eq421}), we have
 \begin{eqnarray*}
 \lambda(H^*,\vec{z})-\lambda(H',\vec{z})&=&z_{l-3}z_{l-2}z_{l}-z_{1}z_{l-1}z_{l}\\
 &=&[(z_{l-2}-z_{l-1})z_{l-3}-(z_1-z_{l-3})z_{l-1}]z_{l}\\
 &>&(z_1-z_{l-3})(z_{l-3}-z_{l-1})z_{l}>0
\end{eqnarray*}
  for $z_{l}>0$. Hence $\lambda(H',\vec{z})<\lambda(H^*,\vec{z})\leq \lambda([l-1]^{(3)})=\lambda(H').$   This contradicts to $\vec{z}$ is an optimal weighting for $H'$. Hence $z_{l-1}=0.$ 
   So we can assume that $H$ is left-compressed. 
Therefore we get a left-compressed extremal $3$-graph $H$ for $m$ and $l-1$. Hence we can assume that $G$ is left -compressed.

Next show that we can assume  $G$ is on $l$ vertices. We will use Lemma \ref{lemma8} below. The proof of Lemma \ref{lemma8} is similar to a proof of a result in \cite{T}. We omit the details.

\begin{lemma}\label{lemma8} Let $m$ and $l$  be  positive integers satisfying ${l-1 \choose 3} \le m \le {l-1 \choose 3} + {l-2 \choose 2}$. Let $G$ be a left-compressed 3-graph on the vertex set $[k]$ with $m$ edges and without containing a clique of order $l-1$  such that $\lambda(G)=\lambda_{(m,l-1)}^{3-}$. Let ${\vec x}$ be an optimal weighting for $G$ with $k$ positive weights. Then Conjecture \ref{conjecture2} holds or
$$\vert [k-1]^{(3)}\setminus E\vert \le k-2.$$
\end{lemma}

\bigskip

Let us continue the proof of Lemma \ref{lemma1}.  Since $G$ is left compressed, then $1(k-1)k\in E$ and $\vert [k-2]^{(2)}\cap E_k \vert \ge 1$.

If $k\ge l+1$, then applying Lemma \ref{lemma8}, we have
\begin{eqnarray*}
m=\vert E\vert &=&\vert E\cap [k-1]^{(3)}\vert +\vert [k-2]^{(2)}\cap E_k \vert +\vert E_{(k-1)k}\vert \\
&\ge & {l \choose 3}-(l-1) +2 \\
&\ge& {l-1 \choose 3} + {l-2 \choose 2}+1,
\end{eqnarray*}
which contradicts to the assumption that $m\le {l-1 \choose 3} + {l-2 \choose 2}$. Recall that $k\ge l$, so we have
$k=l.$ Since ${\vec x}$ has only $l$ positive weights, we can assume that $G$ is on $l$ vertices.
This proves Lemma \ref{lemma1}. \qed

\bigskip

\noindent {\bf Proof of Theorem \ref{theorem2}.}  Since $\lambda_{(m,l-1)}^{3-}$ does not decrease as $m$ increases, it is sufficient to verify Conjecture \ref{conjecture2} for $m = {l-1 \choose 3} + {l-2 \choose 2}$. By Lemma \ref{lemma1}, it is sufficient to verify  $\lambda (G)<\lambda ([l-1]^{(3)})$ for all  left compressed 3-graphs $G$ on $l$ vertices with $m = {l-1 \choose 3} + {l-2 \choose 2}$ edges and without containing the clique $[l-1]^{(3)}$. \qed

\bigskip

Now we describe an Algorithm  to produce all   left compressed 3-graphs $G$ on the vertex set $[l]$ with $m = {l-1 \choose 3} + {l-2 \choose 2}$ edges and without containing the clique $[l-1]^{(3)}$.
Notice that, for a $3$-graph $G$ on $l$ vertices with $m = {l-1 \choose 3} + {l-2 \choose 2}$ edges and without  containing a clique of order $l-1$, we may write $[l]^{(3)}$ as
$$[l]^{(3)}=[l-1]^{(3)}\cup G_1 \cup G_2$$
where $G_1=\{ijl, {\rm where \ } ij\in[l-2]^{(2)}\}$ and $G_2= \{1(l-1)l, 2(l-1)l, 3(l-1)l, \cdots, (l-2)(l-1)l \} $.

Under these assumptions,  $G$ can be obtained from $[l]^{(3)}$ by deletion of a subgraph $H$ with $l-2$ edges since $l-2={l \choose 3}-[{l-1 \choose 3}+{l-2 \choose 2}]$. These $l-2$ edges of $H$ consist of edges from $[l-1]^{(3)}$, $G_1$ (if any), or $G_2$.  Specifically, 3-graph $G$ takes the form of
$$G=([l-1]^{(3)} - E_1) \cup (G_1 - E_2) \cup \{1(l-1)l, 2(l-1)l, 3(l-1)l, \cdots, i(l-1)l \}$$
for some edge set $E_1$ from $[l-1]^{(3)}$, some edge set $E_2$ from $G_1$, and some $i$ where $1 \le i \le l-3$. Observe that $E_{q \setminus j}=\emptyset$ for $ 1\leq q<j\leq i$. By (\ref{enbhd}),  $x_1=x_2= \cdots =x_i $. 

  An $r$-graph $H=(V, E)$ on the vertex set $[l]$ is right-compressed if $i_1i_2i_3\in E$ implies  $j_1j_2j_3\in E$ whenever $j_1\ge i_1$, $j_2\ge i_2$, and $j_3\ge i_3$. Note that $H$ is right-compressed if and only if the complement of $H$ is left-compressed. To generate all possible left compressed $3$-graphs on the vertex set $[l]$  with $m={l-1 \choose 3} + {l-2 \choose 2}$ edges and without containing the clique $[l-1]^{(3)}$, we can generate all right compressed connected $3$-subgraphs $H$ rooted at $(l-2)(l-1)l$ with $l-2$ edges and containing  $(l-3)(l-2)(l-1)$ and then take the complement of each $H$. To do so, we use Algorithm \ref{algorithm1} in the following procedure that is based on Figure 1 (replacing $k$ by $l$).

\begin{algorithm} \label{algorithm1} List all left compressed $3$-graphs on the vertex set $[l]$  with $m={l-1 \choose 3} + {l-2 \choose 2}$ edges and without containing the clique $[l-1]^{(3)}$.

Input: $l\ge 7$ and a Hessian graph with vertex set $[l]^{(3)}$ (replace $k$ by $l$ in  Figure 1).

Output: All right compressed connected $3$-subgraphs $H$ rooted at $(l-2)(l-1)l$ with $l-2$ edges and containing  $(l-3)(l-2)(l-1)$, thus produce all possible left compressed 3-graph $G=[l]^{(3)}-H$ on the vertex set $[l]$ with $m = {l-1 \choose 3} + {l-2 \choose 2}$ edges and without containing the clique $[l-1]^{(3)}$.

Initialization: Set $H=\{(l-2)(l-1)l,(l-3)(l-1)l, (l-3)(l-2)l,(l-3)(l-2)(l-1), (l-4)(l-1)l\}$.

Step 1. For each  direct descendant of  $H$ (an edge $e$ is a  direct descendant of  $H$ if $e\in E(H)$ and $e$ is a direct descendant of an edge in $H$), check whether all its direct ancestors are in $H$. If so, add to $H$. Then record  the new $H$  and record the size of  the new $H$.   Take all distinct new $H$ with size increased by 1 and repeat this process until there are $l-2$ triples in $H$.  Output all distinct $H$ with $l-2$ triples.

Step 2.  Performing $G=[l]^{(3)}-H$ and output $G$.

\end{algorithm}

If $l=6$, then $G=[6]^{(3)}-H$ is the only left-compressed $3$-graph on $[6]$  with ${5 \choose 3}+{4 \choose 2}$ edges without containing a clique of order 5, where $H=\{456, 356, 346, 345\}$.

It is obvious that Algorithm \ref{algorithm1} leads to the following result.

\begin{prop} Algorithm \ref{algorithm1} produces all possible left compressed 3-graphs on $l\ge 7$ vertices with $m = {l-1 \choose 3} + {l-2 \choose 2}$ edges.
\end{prop}






\section{Proof of Corollary  \ref{mbc}}\label{proofmbc}

We need the following result from \cite{SPT}.

\begin{lemma}\label{mainproof2} \cite{SPT} Let $G=(V,E)$ be a left-compressed 3-graph on the vertex set $[l]$ with ${l-1 \choose 3} \le m \le {l-1 \choose 3} + {l-2 \choose 2}$ edges and not containing a clique of order $l-1$.  If $|E_{(l-1)l}|\leq 7$,  then $\lambda (G)< \lambda ([l-1]^{(3)})$.
\end{lemma}


\noindent {\bf Proof of Corollary \ref{mbc}.} Let $6\le l\le 13$. By Theorem \ref{theorem2}, we can assume that
$G$ is a left compressed 3-graph $G$ on vertex set $[l]$ with $m = {l-1 \choose 3}+{l-2 \choose 2}={l \choose 3}-(l-2)$ edges without containing a clique of order $l-1$.  Notice that there are $l-2$ edges in  $G^c$.  Since $G$ is left-compressed and it doesn't contain a clique of order of order $l-1$, then $\{(l-3)(l-2)(l-1), (l-3)(l-2)l, (l-3)(l-1)l, (l-2)(l-1)l\}\subseteq G^c$. Let $i$ be the maximum integer such that $(l-1-i)(l-1)l\in G^c$, then $i\le 9$ since  there are $l-2\le 11$  edges in $G^c$ and $\{(l-3)(l-2)(l-1), (l-3)(l-2)l\}\subset G^c$. Clearly $i\ge 2$.  Note that
\begin{equation}\label{l2i}
\vert E_{(l-1)l}\vert=\vert G^c\cap (\{pql,  pq\in[l-2]^{(2)}\}\cup \{1(l-1)l, 2(l-1)l, 3(l-1)l, \cdots, (l-2)(l-1)l \})\vert=l-2-i.
\end{equation}

If $l\le 11$, then $\vert E_{(l-1)l}\vert \le 9-i\le 7$. By Lemma \ref{mainproof2}, $\lambda (G)< \lambda ([l-1]^{(3)})$ holds.

If $l=12$, then $i\ge 3$. Otherwise, $(l-4)(l-1)l\notin G^c$. Since $G^c$ is right-compressed, there are at most 4 triples ($(l-3)(l-2)(l-1), (l-3)(l-2)l, (l-3)(l-1)l, (l-2)(l-1)l$) in $G^c$. So $i\ge 3$ and  $\vert E_{(l-1)l}\vert\le 7$ (applying (\ref{l2i})). By Lemma \ref{mainproof2}, $\lambda (G)< \lambda ([l-1]^{(3)})$ holds.

If $l=13$, then $i\ge 4$. Otherwise, $(l-5)(l-1)l\notin G^c$. Since $G^c$ is right-compressed, there are at most ${5 \choose 3}=10<11=l-2$ triples  in $G^c$. So $\vert E_{(l-1)l}\vert\le 7$ (applying (\ref{l2i})). By Lemma \ref{mainproof2}, $\lambda (G)< \lambda ([l-1]^{(3)})$ holds.

The proof is completed.\qed

\section{Proof of Theorem \ref{Lemma103}}\label{proofLemma103}

Combining a result from \cite{SPT} and a result from \cite{TPWP}, we have the following lemma.

\begin{lemma} \label{Lemma100} Let $G$  be  a left-compressed $3$-graph on $[l]$ with $m = {l-1\choose 3}+{l-2\choose 2}$ edges. If the first two triples in colex ordering in $G^{c}$  are $(l-3)(l-2)(l-1)$ and $(l-2-i)(l-2)l$ (where $i\ge 1$),  then $ \lambda(G)< \lambda ([l-1]^{(3)})$.
\end{lemma}

\noindent{\bf Proof of Theorem \ref{Lemma103}.  } We apply induction on $j$.
For $j=1$ the assertion is true by Lemma \ref{Lemma100}.  Assume that the assertion is true for $j=s$ (where $s\geq 1)$.  We now consider the case that   $j=s+1$. In this case,  the first triple in colex ordering in $G^c$ is $(l-3-s)(l-2)(l-1)$. Note that $\vert E_{(l-1)l}\vert=i+s+1$.

Let $G'=G\bigcup\{(l-3-s)(l-2)(l-1)\}\backslash \{(i+s+1)(l-1)l\}$, then $\lambda(G') < \lambda ([l-1]^{(3)}) $ by the induction assumption. Next we prove that $\lambda(G)\leq \lambda(G')$.

Let $\vec{x}=(x_{1}, x_{2}, \ldots , x_{l})$  be an optimal weighting for $G$ satisfying $x_1 \ge x_2 \ge \ldots \ge x_l \ge 0$. By Remark \ref{r1}(b), we have $x_1=x_2=x_3=\ldots=x_{i+s+1}$,
\begin{eqnarray}\label{equ1010}
x_{i+s+1}= x_{l-3-s}+\frac{x_{l-2}x_{l}+x_{l-2}x_{l-1}+x_{l-1}x_{l}}{\lambda(E_{(i+s+1)(l-3-s)},\vec{x})},
\end{eqnarray}
and
\begin{eqnarray}\label{equ212}
x_{l-2}=x_l+\frac{(x_{i+s+2}+\ldots+x_{l-4-s})x_{l-1}}{\lambda(E_{(l-2)l},\vec{x})}.
\end{eqnarray}
Note that
\begin{eqnarray}\label{equ10}
\lambda(G',\vec{x})-\lambda(G,\vec{x})=x_{l-3-s}x_{l-2}x_{l-1}-x_{i+s+1}x_{l-1}x_{l}.
\end{eqnarray}

Consider a new weighting for $G'$:  $\vec{y}=(y_1, y_2, \ldots, y_l)$ given by $y_p=x_p$ for $p\neq l-2$, $p\neq l$ and $y_{l-2}=x_{l-2}+\delta $, $y_{l}=x_{l}-\delta$. Note that $\lambda(E'_{(l-2)l},\vec{x}) =\lambda(E_{(l-2)l},\vec{x})$, and $\lambda(E_{(l-2)},\vec{x}) =\lambda(E_{l},\vec{x})$.
Then
\begin{eqnarray*}
\lambda(G',\vec{y})-\lambda(G',\vec{x})&=& \delta(\lambda(E'_{l-2},\vec{x})-\lambda(E'_{l},\vec{x}))- \delta^{2}\lambda(E'_{(l-2)l},\vec{x}) \nonumber\\
&=&\delta[\lambda(E_{l-2},\vec{x})+x_{l-3-s}x_{l-1}-(\lambda(E_{l},\vec{x})-x_{i+s+1}x_{l-1})]-\delta^{2}\lambda(E'_{(l-2)l},\vec{x})\nonumber \\
&=&\delta(x_{l-3-s}x_{l-1}+x_{i+s+1}x_{l-1})-\delta^{2}\lambda(E_{(l-2)l},\vec{x}).
\end{eqnarray*}
Let $$\delta=\frac{ x_{l-3-s}x_{l-1}+x_{i+s+1}x_{l-1}}{2\lambda(E_{(l-2)l},\vec{x})}.$$
Clearly,
$$\delta \le \frac{x_{l-3-s}x_{l-1}+x_{i+s+1}x_{l-1}}{2(i+s+1)x_{1}}\leq\frac{x_{l-1}}{i+s+1}.$$
By Remark \ref{r1}(b),  if $i=s+1$, then
\begin{eqnarray}
x_{l-1}=x_{l},
\end{eqnarray}
and $\delta\leq\frac{x_{l-1}}{i+s+1}<x_l$.

If $i\geq s+2$, then
\begin{eqnarray*}
x_{l-2}=x_{l-1}+\frac{(x_{i+s+2}+\ldots+x_{l-3-i})x_{l}}{\lambda(E_{(l-2)(l-1)},\vec{x}}=x_{l-1}+\frac{(x_{i+s+2}+\ldots+x_{l-3-i})x_{l}}{x_{1}+\ldots+x_{l-4-s}}
\leq x_{l-1}+x_{l},
\end{eqnarray*}
and
\begin{eqnarray*}
x_{l-1}=x_{l}+\frac{(x_{l-2-i}+\ldots+x_{l-4-s})x_{l-2}}{(i+s+1)x_1} \leq x_{l}+\frac{(i-s-1)x_{l-2}}{i+s+1}\leq x_{l}+\frac{i-s-1}{i+s+1}x_{l-1}+\frac{i-s-1}{i+s+1}x_{l}.\nonumber\\
\end{eqnarray*}
So $\frac{x_{l-1}}{i+s+1}\leq \frac{2ix_{l}}{(2s+2)(i+s+1)}=\frac{ix_{l}}{(s+1)(i+s+1)}<x_{l}$. Recall that  $\delta\leq\frac{x_{l-1}}{i+s+1}$. Therefore, $\delta\leq x_{l}$. Hence $\vec{y}=(y_1, y_2, \ldots, y_l)$ is also a legal weighting for $G'$. So we have
\begin{eqnarray}\label{equ1014}
\lambda(G',\vec{y})-\lambda(G',\vec{x})=\frac{(x_{l-3-s}+x_{i+s+1})^{2}x_{l-1}^{2}}{4\lambda(E_{(l-2)l},\vec{x})}.
\end{eqnarray}

Let $\vec{y'}=(y'_1, y'_2, \ldots, y'_l)$ given by $y'_q= y_q$ for $q\neq i+s+1$, $q\neq l-3-s$ and $y'_{i+s+1}=y_{i+s+1}-\eta$, $y'_{l-3-s}=y_{l-3-s}+\eta$. Note that $\lambda(E'_{(i+s+1)(l-3-s)},\vec{y})=\lambda(E_{(i+s+1)(l-3-s)},\vec{x})$. Then
\begin{eqnarray*}
\lambda(G',\vec{y'})-\lambda(G',\vec{y})&=& \eta (\lambda(E'_{l-3-s},\vec{y})-\lambda(E'_{i+s+1},\vec{y}))-\eta ^{2}\lambda(E'_{(l-3-s)(i+s+1)},\vec{ y}) \nonumber\\
&=& \eta (y_{i+s+1}- y_{l-3-s})\lambda(E'_{(i+s+1)(l-3-s)},\vec{y})-\eta^{2}\lambda(E'_{(i+s+1)(l-3-s)},\vec{y}).
\end{eqnarray*}
Let $$ \eta=\frac{y_{i+s+1}-y_{l-3-s}}{2}=\frac{x_{i+s+1}-x_{l-3-s}}{2}.$$
Clearly  $\vec{y'}=(y'_1, y'_2, \ldots, y'_l)$ is also a legal weighting for $G'$. Applying (\ref{equ1010}), we have
\begin{eqnarray}\label{equ1015}
\lambda(G',\vec{y'})-\lambda(G',\vec{y})= \frac{(x_{l-2}x_{l}+x_{l-2}x_{l-1}+x_{l-1}x_{l})^{2}}{4\lambda( E_{(i+s+1)(l-3-s)},\vec{x})}.
\end{eqnarray}
Combing  (\ref{equ1010}), (\ref{equ212}), (\ref{equ10}), (\ref{equ1014}) and (\ref{equ1015}), we have
\begin{eqnarray*}
\lambda(G',\vec{y'})-\lambda(G,\vec{x})&=&\frac{(x_{l-2}x_{l}+x_{l-2}x_{l-1}+x_{l-1}x_{l})^{2}}{4\lambda( E_{(i+s+1)(l-3-s)},\vec{x})}
+\frac{(x_{l-3-s}+x_{i+s+1})^{2}x_{l-1}^{2}}{4\lambda(E_{(l-2)l},\vec{x})}\\
&&+x_{l-3-s}x_{l-2}x_{l-1}-x_{i+s+1}x_{l-1}x_{l}\nonumber
\\&=&\frac{(x_{l-2}x_{l}+x_{l-2}x_{l-1}+x_{l-1}x_{l})^{2}}{4\lambda(E_{(i+s+1)(l-3-s)},\vec{x})}+\frac{(x_{l-3-s}+x_{i+s+1})^{2}x_{l-1}^{2}}{4\lambda(E_{(l-2)l},\vec{x})}
\\&+&x_{l-3-s}x_{l-1}(x_{l-2}-x_{l})-(x_{i+s+1}-x_{l-3-s})x_{l-1}x_{l}\\&=&\frac{(x_{l-2}x_{l}+x_{l-2}x_{l-1}+x_{l-1}x_{l})^{2}}{4\lambda(E_{(i+s+1)(l-3-s)},\vec{x})}
+\frac{(x_{l-3-s}+x_{i+s+1})^{2}x_{l-1}^{2}}{4\lambda(E_{(l-2)l},\vec{x})}
\\&+&\frac{x_{l-3-s}x_{l-1}^{2}(x_{i+s+2}+\ldots+x_{l-4-s})}{\lambda(E_{(l-2)l},\vec{x})}
-\frac{x_{l-2}x_{l-1}^{2}x_{l}+x_{l-2}x_{l-1}x_{l}^{2}+x_{l-1}^{2}x_{l}^{2}}{\lambda(E_{(i+s+1)(l-3-s)},\vec{x})}
\\&\geq&\frac{x_{l-2}x_{l-1}x_{l}^2+x_{l-1}^{2}x_l^2-x_{l-2}x_{l-1}^{2}x_{l}-x_{l-2}x_{l-1}x_{l}^{2}-x_{l-1}^{2}x_{l}^{2}}{\lambda(E_{(i+s+1)(l-3-s),\vec{x}})}
\\&+&\frac{x_{l-3-s}x_{i+s+1}x_{l-1}^{2}+x_{l-3-s}x_{l-1}^{2}(x_{l-4-s}+\ldots+x_{i+s+2})}{\lambda(E_{(l-2)l},\vec{x})}
\\&\geq&\frac{x_{l-3-s}x_{i+s+1}x_{l-1}^{2}+x_{l-3-s}x_{l-1}^{2}(x_{l-4-s}+\ldots+x_{i+s+2})-x_{l-2}x_{l-1}^{2}x_{l}}{\lambda(E_{(l-2)l},\vec{x})}\geq 0.
\end{eqnarray*}
So $\lambda(G)=\lambda(G,\vec{x})\leq \lambda(G',\vec{y'})\leq \lambda(G')< \lambda([l-1]^{(3)})$.

\qed


%

\end{document}